\documentclass{article}
\usepackage[all]{xy}
\usepackage{enumerate}
\usepackage{amsmath,amsfonts,amsthm,amssymb}
\newtheorem{theorem}{Theorem}[section]
\newtheorem{example}[theorem]{Example}
\newtheorem{definition}[theorem]{Definition}
\newtheorem{proposition}[theorem]{Proposition}

\newtheorem{corollary}[theorem]{Corollary}
\newtheorem{remark}[theorem]{Remark}
\renewcommand{\Im}[1]{\textnormal{Im}(#1)}
\newcommand{\im}[1]{\textnormal{im}(#1)}
\newcommand{\C}{\mathbb{C}}

\newcommand{\I}{\mathbb{I}}

\newcommand{\sq}{\sqsubset}
\newcommand{\E}{\mathcal{E}}
\newcommand{\M}{\mathcal{M}}
\newcommand{\fac}{(\E,\M)}
\newcommand{\facf}[1]{(\E^#1,\M^#1)}
\newcommand{\cod}{\textnormal{cod}}
\newcommand{\two}{\textbf{2}}
\begin{document}
\title{A note on images of cover relations}
\author{J. R. A. Gray}
\maketitle
\begin{abstract}
For a category $\C$, a
small category $\I$, and a pre-cover relation $\sq$ on $\C$
we prove, under certain completeness assumptions on $\C$, that
a morphism $g: B\to C$ in the functor category $\C^\I$ admits an image with
respect to the pre-cover relation on $\C^\I$ induced by $\sq$ as soon as each
component of $g$ admits an image with respect to $\sq$.
 We then apply this to show that if a pointed category $\C$ is: (i)
 algebraically cartesian closed; (ii) exact protomodular and action accessible; or
 (iii) admits normalizers, then the same is true of each functor category $\C^\I$
with $\I$ finite. In addition, our results give explicit constructions of
 images in
functor categories using limits and images in the underlying category.
 In particular, they can be used to give explicit constructions of
both centralizers and normalizers in functor categories using limits and
 centralizers or normalizers (respectively) in the underlying category.
\end{abstract}
\section*{Introduction}
Perhaps the most natural way to extend the definition of commuting elements of
a group to homomorphisms into a group, is to say that a pair of group
homomorphisms $f:A\to C$ and $g:B\to C$ commute if each element in the image
of $f$ commutes with each element in the image of $g$. Moreover, one can show,
directly or as a special case of S.{} Mac Lane's
characterization of bifunctors 
(see 
\cite{MAC_LANE:1997}), that this condition is equivalent to the existence
of a morphism $\varphi: A\times B\to C$ making the diagram
\begin{equation}
\label{diagram:commutes}
\vcenter{
\xymatrix{
A \ar[r]^-{\langle 1,0\rangle}\ar@/_3ex/[dr]_{f} & A\times B \ar[d]^{\varphi}&
B\ar[l]_-{\langle 0,1\rangle}\ar@/^3ex/[dl]^{g}\\
& C, &
}
}
\end{equation}
in which $\langle 1,0\rangle$ and $\langle 0,1\rangle$ are the homomorphisms
defined by $\langle 1,0 \rangle(a)=(a,1)$ and $\langle 0,1\rangle(b)=(1,b)$
respectively, commute. This last formulation was used by S.{} Huq
in \cite{HUQ:1968} to study commutativity and other closely related notions
in a categorical context, close to the more recent semi-abelian context
introduced by G.{} Janelidze, L.{} M\'arki, and T.{} Tholen in
\cite{JANELIDZE_MARKI_THOLEN:2002}.
Later, Z.{} Janelidze \cite{JANELIDZE_Z:2009}
introduced and studied relations on the morphisms of a category
called cover relations and of which the \emph{commutes relation}
described above is an example. It was shown that other similar kinds of
cover relations arise from special kinds of monoidal structures which are
called monoidal sum structures. In addition, it was shown that cover relations
also arise from factorization systems on categories, from which some of
terminology and notation is derived.

Let us recall briefly what a cover relation is and how both these kinds of
cover relations arise.
A pre-cover relation $\sq$ is a relation on the class of
morphisms of
a category $\C$ such that if $f\sq g$ then the codomain of $f$ and $g$ are
equal. A pre-cover relation $\sq$ is a cover relation if it satisfies:
(i) if $f \sq g$ and the composite $hf$ is defined, then $hf \sq hg$;
(ii) if $f\sq g$ and the composite $fe$ is defined, then $fe \sq g$.
Let $(\C,\otimes,I,\alpha,\rho,\lambda)$ be a monoidal category, such that
$I$ is an initial object in $\C$ and for each $A$ and $B$ the morphisms
\[
\xymatrix{
A\otimes I \ar[r]^-{1_A\otimes !_B} &  A\otimes B  &
I\otimes B, \ar[l]_-{!_A\otimes 1_B}
}
\] 
where $!_A$ and $!_B$ are the unique morphism from $I$ to $A$ and $B$
respectively, are jointly epimorphic. The induced relation $\sq$ is defined
by requiring, for a pair of morphisms $f:A\to C$ and $g:B\to D$,
that $f\sq g$ whenever $C=D$ and there exists a morphism
$\varphi : A\otimes B \to C$ making the diagram
\[
\xymatrix{
A \ar[r]^-{\rho_A}\ar@/_2ex/[drr]_{f} & A\otimes I \ar[r]^-{1_A\otimes !_B} &
 A\otimes B\ar[d]^{\varphi}  & I\otimes B, \ar[l]_-{!_A\otimes 1_B} &
B\ar[l]_-{\lambda_B}\ar@/^2ex/[dll]^-{g}\\
&& C &&
}
\]
commute. On the other hand given a factorization system $\fac$ on a
category $\C$ the induced relation $\sq$ on the morphisms of $\C$ is defined
by requiring $f \sq g$ whenever $f$ and $g$ have the same codomain and
if $g = mv$ with $m$ in $\M$, then there exists $u$ such that $f=mu$. 
One of the aims of \cite{JANELIDZE_Z:2009} was to show that under
suitable conditions, both factorization systems and monoidal structures can be
recovered from their induced cover relations. 

Let us also recall that,
the image of a morphism $g:B\to C$ with respect to a cover relation $\sq$
on $\C$
can be defined
as the terminal
object in the full subcategory of $(\C\downarrow C)$ with objects $(A,f)$
such that $f\sq g$. A simple observation in recovering a factorization system
$\fac$ on a category
$\C$ from its induced cover relation $\sq$
is to note that
if the class $\M$ of a factorization system $\fac$ consists of monomorphisms,
then it is the class of all morphisms
$f:A\to C$ such that $(A,f)$ is the image of some morphism $g$ with respect
to $\sq$. In this case it turns out that if $(A,f)$ is the image of
a morphism $g : B\to C$
with respect to induced cover relation, then $f$ is the image of $g$ with
respect to the
factorization system. Given a factorization system $\fac$ on a category $\C$,
such that $\M$ consists of monomorphisms, it is easy to observe that for each
category $\I$ the induced factorization system $\facf{\I}$ on the
functor category $\C^\I$, defined componentwise, has $\M^\I$ consisting of
monomorphisms. This means that the induced
cover relation $\sq^\I$ on $\C^\I$, defined componentwise, admits images since
it is also the cover relation induced by $\facf{\I}$, and furthermore, its
images are \emph{computed componentwise}. However, it is not the case that
if the underlying cover relation admits images, then
images necessarily exist for the induced cover relations on each functor
category
see Example
\ref{example:images_don't lift} below, nor that when they do exist they
are computed componentwise. Our
aims here are: (a) to prove that if $\C$ and $\I$ are categories and $\sq$ is
a cover relation (more generally a \emph{pre-cover relation} satisfying
Condition \ref{def:1} (ii) below) on $\C$
admitting images, then, under certain completeness
assumptions on $\C$, the induced cover relation (pre-cover relation)
$\sq^{\I}$ on $\C^\I$ 
admits
images; (b) to apply this result to prove that several categorical
algebraic conditions lift from a category to its functor categories. In
particular, we show that for a \emph{semi-abelian category} $\C$
\emph{action accessibility}, \emph{algebraically-cartesian closedness},  and
the existence of normalizers, \emph{lift} to functor categories  $\C^{\I}$
with $\I$ finite. In fact our results hold more generally see Corollary
\ref{main:corollary} for a precise formulation.
\section{Preliminaries}
In this section we recall the basic background on (pre-)cover relations that
we will
need in the next section. In addition, we give examples of (pre-)cover
relations
(one of which is new), to which we will apply our main theorem of the
following section, to obtain the above mentioned results showing that
certain categorical algebraic conditions lift to functor categories.
\begin{definition}[\cite{JANELIDZE_Z:2009}]\label{def:1}
A pre-cover relation $\sq$ is a relation on class of morphisms of a category
$\C$, such that if $f\sq g$ then the codomain of $f$ and $g$ are equal.
A pre-cover relation $\sq$ is a cover relation if it satisfies:
\begin{enumerate}[(i)]
\item if $f \sq g$ and the composite $hf$ is defined, then $hf \sq hg$;
\item if $f\sq g$ and the composite $fe$ is defined, then $fe \sq g$.
\end{enumerate}
\end{definition}
\begin{definition}[\cite{JANELIDZE_Z:2009}]
Let $\sq$ be a pre-cover relation on a category $\C$.
The image of a morphism $g:B\to C$ is terminal object in the full subcategory
 of $(\C\downarrow C)$ with objects $(A,f)$ such that $f\sq g$.
 We will denote the image of $g$ by $(\Im{g},\im{g})$.
\end{definition}
One easily observes:
\begin{proposition}[\cite{JANELIDZE_Z:2009}]
\label{proposition:images_are_monos_for_right_cover_relations}
Let $\sq$ be a pre-cover relation on $\C$ satisfying
Condition \ref{def:1} (ii)
and let $g:B\to C$ be a morphism in $\C$.
If $g$ admits an image, then the morphism $\im{g}$ is a monomorphism.
\end{proposition}
Let us recall the necessary background in order to give the two examples of
cover relations mentioned above. For pointed category $\C$ we
write $0$ for the zero object as well as for each zero morphism between
each pair of objects. For objects $A$ and $B$ we write
$\pi_1: A\times B\to A$ and $\pi_2 : A\times B\to B$ for the first and second
product projections (whenever they exist), and for a pair of morphisms
$f:W\to A$ and $g:W\to B$ we write $\langle f,g\rangle : W\to A\times B$
for the unique morphism with 
$\pi_1\langle f,g\rangle = f$ and $\pi_2\langle f,g\rangle  =g$.
Recall that a pointed category with finite limits is called unital
\cite{BORCEUX_BOURN:2004}, if
for objects $A$ and $B$ the morphisms $\langle 1,0\rangle: A\to A\times B$ and
$\langle 0,1\rangle: B \to A\times B$ are jointly extremal-epimorphic. A pair
of morphisms $f:A\to C$ and $g: B\to C$ are said to (Huq) commute if there
exists a (necessarily unique) morphism $\varphi : A\times B\to C$ making the
diagram \eqref{diagram:commutes} commute.
\begin{example}\label{example:commutes}
 The commutes relation on a unital category is a cover relation.  (In fact it
 is the cover relation induced by the cartesian monoidal structure on $\C$).
 The image of a morphism $g$ with respect to this cover relation is called the
centralizer of $g$.
\end{example}
\begin{example}\label{example:normalizes}
Let $\C$ be a pointed category. For morphisms $f:A\to X$ and $g:B\to Y$ in
 $\C$ we say that $f$
normalizes $g$ if $X=Y$,
and there exists a morphism $u:A\to C$,
 a normal monomorphism $v:B\to C$, and monomorphism $h:C\to X$
 making the diagram
 \[
  \xymatrix{
   A\ar[dr]_{f} \ar[r]^{u}& C \ar[d]^{h} &B\ar[dl]^{g}\ar[l]_{v}
\\
    & X &   }
 \]
 commute.
 Note that in this case, it follows that $g$ is a monomorphism and
 that $h$ normalizes
 $g$ too since $h = h 1_C$ and $g=hv$. The normalizes relation is a pre-cover
relation satisfying Condition \ref{def:1} (ii)
and the image of a monomorphism $g:B\to C$ is the normalizer of $g$
 (in the sense of \cite{GRAY:2013b}).
 Indeed, if $g:B\to X$ is a monomorphism and $(\Im{g},\im{g})$ is the image of $g$,
  then there exists a morphism $u : \Im{g} \to C$, a monomorphism $h : C\to X$ and a normal
  monomorphism $v:B\to C$ such that $\im{g} = hu$ and $g=hv$. But, since $h$
  normalizes $g$ it follows that $u : (\Im{g},\im{g}) \to (C,h)$ being a monomorphism with
  domain a terminal object is an isomorphism.
\end{example}
\begin{remark}
 When $\C$ is \emph{ideal determined} \cite{JANELIDZE_MARKI_THOLEN_URSINI:2010}
 the pre-cover relation in the previous
 example can be enlarged to
 become a cover relation as follows.
 For $f$ and $g$ as above instead of requiring the morphism $v$ (required to
 exits) to be a normal monomorphism we require that its (regular) image must be a
 normal monomorphism.
 The idea of considering this cover relation arose in discussions with
 Z.\ Janelidze and led to considering the above pre-cover relation.
\end{remark}
\section{The main results}
Let $\sq$ be a pre-cover relation on a category $\C$. For a category $\I$, the
pre-cover relation $\sq$ induces a pre-cover relation $\sq^\I$ on $\C^\I$,
which is defined
componentwise. Note that if $\sq$ satisfies Condition \ref{def:1} (i) or (ii)
then so does $\sq^\I$. We have:
\begin{theorem}
Let $\I$ be a category and let $\C$ be a category with pullbacks and
with wide pullbacks of families of monomorphisms indexed by the morphisms of
$\I$. Let $\sq$ be a pre-cover relation on $\C$ satisfying Condition \ref{def:1} (ii).
A morphism $g:B\to C$ in the functor category $\C^\I$, admits an image
with respect to $\sq^\I$ on $\C^\I$ if each component of $g$ admits an image
with respect to $\sq$. Moreover, when this is the case, for each $X$ in $\I$
  the object $(\Im{g}(X),\im{g}_X)$ is the product in
 the comma category $(\C\downarrow C(X))$
 of a certain family $(W_i,w_i)_{i\in \mathcal I}$ where $\mathcal I$ is the
 collection of all morphisms with domain $X$.
 This family consists of monomorphisms $w_i: W_i\to C(X)$
 obtained for each $i:X\to Y$ in $\mathcal I$ by pulling back
 $\im{g_{Y}}$ along $C(i)$,
 as displayed in the lower square of the diagram
\begin{equation}
\label{diagram:image_in_functor_category}
\vcenter{
\xymatrix{
\Im{g}(X) \ar[drr]^{\im{g}_X} \ar[dr]_{v_i} & &\\
& W_i \ar[r]^{w_i} \ar[d]_{\tilde{i}}& C(X)\ar[d]^{C(i)}\\
& \Im{g_{Y}} \ar[r]_{\im{g_{Y}}} & C(Y).
}
}
\end{equation}
\end{theorem}
\begin{proof}
Let $g: B\to C$ be morphism in $\C^\I$.  We begin by showing that the
 above mentioned construction produces a
morphism $f:A\to C$ which we then show to be the image of $g$. For each
$i:X\to Y$ in $\I$, let $w_i: W_i\to C(X)$ be the preimage of $\im{g_{Y}}$
along $C(i)$ as displayed in \eqref{diagram:image_in_functor_category}.
For each object $X$ in $\I$, let $v_i : A(X)\to W_i$ be the $i$th projection
of the wide
pullback of all $w_i$ where $i$ is a morphism with domain $X$, and let
$f_X  = w_{1_X} v_{1_X}$.  Now let $i: X\to Y$ be a morphism in $\I$.
 It is easy to check that for each morphism $j:Y\to Z$ there exists a unique
morphism $\bar i: W_{ji} \to W_{j}$ making the right hand square in the
diagram
\[
\xymatrix{
A(X)\ar@{-->}[d]_{A(i)} \ar[r]^{v_{ji}} & 
W_{ji} \ar[r]^{w_{ji}} \ar[d]_{\bar i}& C(X) \ar[d]^{C(i)}\\
A(Y) \ar[r]^{v_{j}} & W_{j} \ar[r]_{w_{j}} & C(Y) 
}
\]
commute (in fact making it a pullback). 
It follows that there exists a unique morphism $A(i):A(X)\to A(Y)$ such that
the left hand diagram above commutes for each such $j$. These assignments
make $A$ an object and $f: A \to C$ a morphism in $\C^\I$.  Since, by
definition, for each $X$ in $\I$ the diagram 
\[
\xymatrix{
A(X) \ar[r]_{v_{1_X}} \ar@/^3ex/[rr]^{f_{X}} &
W_{1_X} \ar[r]_{w_{1_X}} \ar[d]_{\tilde{1}_X}& C(X)\ar[d]^{C(1_X)}\\
& \Im{g_{X}} \ar[r]_{\im{g_{X}}} & C(X).
}
\]
commutes, we see that $f_X = \im{g_X} \tilde 1_X v_{1_X}$ and
so $f_X \sq g_X$ and $f \sq g$.  Let $f' : A' \to C$ be a morphism in
$\C ^\I$ such that $f' \sq g$, we need to show that there exists a unique
morphism $u : A' \to A$ such that $fu = f'$. Since $f'_{Y} \sq g_{Y}$,
for each morphism $i:X\to Y$, there exists a unique morphism
$\underline u_{Y}: A'(Y) \to \Im{g_{Y}}$ such that the solid arrows in the 
diagram
\[
\xymatrix{
A'(X) \ar@{-->}[r]_{v'_i} \ar@/^3ex/[rr]^{f'_X}\ar[d]_{A'(i)} & 
 W_i \ar[d]_{\tilde{i}} \ar[r]_{w_i} & C(X) \ar[d]^{C(i)}\\
A'(Y) \ar[r]^{\underline u_{Y}} \ar@/_3ex/[rr]_{f'_{Y}} 
& \Im{g_{Y}} \ar[r]^{\im{g_{Y}}} & C(Y) 
}
\]
commute, and hence there exists a unique morphism $v'_i : A'(X) \to W_i$
making the entire diagram commute. It now follows that there exists a unique
morphism $u_{X} : A'(X) \to A(X)$ such that $v_i u_{X} = v'_i$ for each $i$
with domain $X$ and hence such that $f_Xu_X=f'_X$. Noting that the
components of $f$ are monomorphisms it follows that the morphisms $u_X$  are
components of a (unique) natural transformation $u:A'\to A$ with $fu=f'$, as
required.
\end{proof}
\begin{remark}
 Suppose that $\I$ and $\C$ are categories and 
 let $U: \I_0\to \I$ be the functor including the objects of $\I$
 as a discrete category in $\I$. Recall that the induced functor
 $\C^U : \C^{\I}\to \C^{\I_0}$ has right adjoint given by taking right
 Kahn extensions, and these Kahn extensions are computed pointwise when
 $\C$ admits certain products (for instance when $\C$ admits products of
 families whose indexing set is bounded by the morphisms of $\I$). 
 Recall also that if $\C$ has pullbacks and the functor $\C^U$ has a right
 adjoint, then for each $C$ in $\C^{\I}$ the induced functor
 $\C^U_{C}:(\C^\I \downarrow C)\to (\C^{\I_0}\downarrow CU)$ also has a
 right adjoint. Note that each functor $\C^U_C$ also has a right adjoint
 (although the functor $\C^U$ may not)
 if $\C$ admits wide pullbacks of families whose indexing set is bounded
 by the morphisms of $\I$. Now suppose that $g:B\to C$ is a morphism in
 $C^\I$ such that $g_X$ admits an image for each $X$ in $\I$,
 and $\C^U_{C}$ has right adjoint $R$. One easily checks that the
 image of $gU$ in $\C^{\I_0}$ is computed componentwise and that
 $(\Im{g},\im{g}) = R(\Im{gU},\im{gU})$.
\end{remark}
Recall that a category $\C$ is semi-abelian
\cite{JANELIDZE_MARKI_THOLEN:2002} in the sense of G.{} Janelidze,
L.{} M\'arki, and W.{} Tholen if it is pointed, (Barr)-exact \cite{BARR:1971},
Bourn-protomodular \cite{BOURN:1991} and has binary coproducts.
\begin{corollary}\label{main:corollary}
Let $\C$ be a finitely complete category and let $\I$ be a finite category.
 If $\C$ is
 \begin{itemize} 
  \item  unital and algebraically cartesian closed
 \cite{BOURN_GRAY:2012} (see also \cite{GRAY:2010b,GRAY:2012b} where this
 notion is first considered but unnamed), or
 \item semi-abelian (more generally pointed exact
protomodular) and action accessible \cite{BOURN_JANELIDZE:2009}, or
   \item admits
 normalizers \cite{GRAY:2013b},
 \end{itemize}
 then the same is true of the functor category
$\C^\I$.
\end{corollary}
\begin{proof}
Recall that: (i) a unital category is algebraically cartesian closed if and
only if it admits centralizers (in the sense above), 
 Proposition 1.2 of \cite{BOURN_GRAY:2012}; (ii) a pointed exact
protomodular category is action accessible if and only if for each normal
monomorphism $n:S\to C$ the normalizer of
$\langle n,n\rangle : S\to C\times C$ exists, Theorem 3.1 of \cite{GRAY:2015a},
 (see also \cite{BOURN_GRAY:2015}).
 The claim now
follows from the previous theorem applied to the pre-cover relations in
Examples \ref{example:commutes} and
\ref{example:normalizes}. Just note that the conditions of being pointed,
unital, exact, and protomodular
easily lift to functor categories.
\end{proof}
We end the paper by giving a simple example showing that images don't always
lift to functor categories.
\begin{example}
\label{example:images_don't lift}
Let $\C$ be the poset with underlying set with distinct elements
$\{A, A', B, B', C, C'\}$ and
with partial order generated from
 \[\{(A,A'), (A,C), (B,B'), (B,C), (A',C'), (B',A'), (C,C')\},\] considered as
a category. For convenience, we introduce labels for some of the morphisms
as shown in the diagram
\[
\xymatrix{
A \ar[d]_{\alpha}\ar[r]^{f} & C \ar[d]^{\gamma} & B\ar[l]_{g} \ar[d]^{\beta}\\
A' \ar[r]^{f'} & C' & B'.\ar[l]_{g'}\ar@/^3ex/[ll]^{u}
}
\]
Now let
$\sq = \{(f'\sigma,g')\,|\, \sigma \in \C_1, \cod(\sigma) = A'\} 
\cup \{(\theta,\phi) \in \C_1^2\,|\, \cod(\theta) = \cod(\phi), \phi\neq g'\}$
 where $\C_1$
 is the set of morphisms of $\C$, and $\cod$ is the \emph{codomain map}.
It is easy to check that $\sq$ is a cover
relation. The main point is to note that if $(\theta,\phi)$ is in $\sq$,
$h\phi = g'$ and $h\neq 1_{C'}$,
then either $h=f'$ and $\phi=u$, or $h=g'$ and $\phi=1_{B'}$ and hence in
either case $h\theta = f'\sigma$ for some $\sigma$.  It straightforward to check that
the image of a morphism $\phi$ is $(\cod{\phi},1_{\cod{\phi}})$, unless
$\phi=g'$ in which case its image is $(A',f')$.
However, the morphism
$(g,g') : (B,B',\beta) \to (C,C',\gamma)$ in $\C^{\two}$, the category of
morphisms of $\C$, has no image with respect to $\sq^{\two}$.  To see why,
 just note
 that the full subcategory of $(\C^{\two}\downarrow (C,C',\gamma))$, with
objects all objects $((X,X',\chi),(p,p'))$ in 
$(\C^{\two}\downarrow (C,C',\gamma))$
 such that $((p,p'),(g,g')) \in \sq^{\two}$, considered as poset, has
 two maximal elements
 $((A,A',\alpha),(f,f'))$ and
 $((B,A',u\beta),(g,f'))$ (and hence no largest element).
\end{example}
\providecommand{\bysame}{\leavevmode\hbox to3em{\hrulefill}\thinspace}
\providecommand{\MR}{\relax\ifhmode\unskip\space\fi MR }
\providecommand{\MRhref}[2]{%
  \href{http://www.ams.org/mathscinet-getitem?mr=#1}{#2}
}
\providecommand{\href}[2]{#2}


\begin{thebibliography}{10}

\bibitem{BARR:1971}
M.~Barr, \emph{Exact categories}, in: Lecture Notes in Mathematics
  \textbf{236}, 1--120, 1971.

\bibitem{BORCEUX_BOURN:2004}
F.~Borceux and D.~Bourn, \emph{Mal'cev, protomodular, homological and
  semi-abelian categories}, Kluwer Academic Publishers, 2004.

\bibitem{BOURN:1991}
D.~Bourn, \emph{Normalization equivalence, kernel equivalence and affine
  categories}, Lecture Notes in Mathematics, Category theory ({C}omo, 1990)
  \textbf{1488~}, Springer, Berlin, 43--62, 1991.

\bibitem{BOURN_GRAY:2012}
D.~Bourn and J.~R.~A. Gray, \emph{Aspects of algebraic exponentiation},
  Bulletin of the Belgian Mathematical Society \textbf{19}(5), 821--844, 2012.

\bibitem{BOURN_GRAY:2015}
D.~Bourn and J.~R.~A. Gray, \emph{Normalizers and split extensions}, Applied
  Categorical Structures \textbf{23}(6), 753--776, 2015.

\bibitem{BOURN_JANELIDZE:2009}
D.~Bourn and G.~Janelidze, \emph{Centralizers in action accessible categories},
  Cahiers de Topologie et G\'eom\'etrie Diff\'erentielles Cat\'egoriques
  \textbf{50}(3), 211--232, 2009.

\bibitem{GRAY:2010b}
J.~R.~A. Gray, \emph{Algebraic exponentiation and internal homology in general
  categories}, Ph.D. thesis, University of Cape Town, 2010.

\bibitem{GRAY:2012b}
J.~R.~A. Gray, \emph{Algebraic exponentiation in general categories}, Applied
  Categorical Structures \textbf{20}(6), 543--567, 2012.

\bibitem{GRAY:2013b}
J.~R.~A. Gray, \emph{Normalizers, centralizers and action representability in
  semi-abelian categories}, Applied Categorical Structures \textbf{22}(5-6),
  981--1007, 2014.

\bibitem{GRAY:2015a}
J.~R.~A. Gray, \emph{Normalizers, centralizers and action accessibility},
  Theory and Applications of Categories \textbf{30}(12), 410--432, 2015.

\bibitem{HUQ:1968}
S.~A. Huq, \emph{Commutator, nilpotency and solvability in categories},
  Quarterly Journal of Mathematics \textbf{19}(1), 363--389, 1968.

\bibitem{JANELIDZE_MARKI_THOLEN:2002}
G.~Janelidze, L.~M\'arki, and W.~Tholen, \emph{Semi-abelian categories},
  Journal of Pure and Applied Algebra \textbf{168}, 367--386, 2002.

\bibitem{JANELIDZE_MARKI_THOLEN_URSINI:2010}
G.~Janelidze, L.~M\'arki, W.~Tholen, and A.~Ursini, \emph{Ideal determined
  categories}, Cahiers de topologie et geom\'etrie diff\'erentielle
  cat\'egoriques \textbf{51}, 115--127, 2010.

\bibitem{JANELIDZE_Z:2009}
Z.~Janelidze, \emph{Cover relations on categories}, Applied Categorical
  Structures \textbf{17}(4), 351--371, 2009.

\bibitem{MAC_LANE:1997}
S.~{Mac Lane}, \emph{Categories for the working mathematician}, 2nd edition
  ed., Springer Science, 1997.

\end{thebibliography}
\end{document}